\documentclass[12pt]{amsart}
\usepackage{amsthm}
\usepackage{amsmath,amssymb,amsfonts,amsthm,color,latexsym}
\usepackage{xcolor}
\usepackage{float}

\newcommand{\C}{\mathbb{C}}
\newcommand{\F}{\mathcal{F}}

\newcommand{\bpartial}{\bar{\partial}}
\newcommand{\dif}{\textsf{Diff}}

\newcommand{\sing}{\textsf{Sing}}

\newcommand{\cod}{\textsf{cod}}
\newtheorem{theorem}{Theorem}[section]

\newtheorem{maintheorem}{Theorem}

\usepackage[pagebackref]{hyperref}

\newtheorem{exemplo}{Example}

\newtheorem{lemma}[theorem]{Lemma}

\newtheorem{definition}{Definition}[section]

\newtheorem{remark}{Remark}[section]

\date{}

\title[A Morse-Bott normal form for Levi-flat hypersurfaces]{A Morse-Bott normal form for real analytic Levi-flat hypersurfaces}
\author{Arturo Fern\'andez-P\'erez}
\address[A. Fern\'andez-P\'erez]{Departamento de Matem\'atica, Universidade Federal de Minas Gerais, UFMG}
\curraddr{Av. Pres. Ant\^onio Carlos 6627, 31270-901, Belo Horizonte-MG, Brazil.}
\email{fernandez@ufmg.br}
\author{Gustavo Marra}
\address[Gustavo Marra]{Institutos de Ci\^encias Puras e Aplicadas, Universidade Federal de Itajub\'a, UNIFEI}
\curraddr{Rua Irm\~a Ivone Drumond 200, 35903-087, Itabira-MG, Brazil.}
\email{marra@unifei.edu.br}

\thanks{The first author acknowledges support from CNPq Projeto Universal 408687/2023-1 ``Geometria das Equa\c{c}\~oes Diferenciais Alg\'ebricas" and CNPq-Brazil PQ-306011/2023-9}

\begin{document}
\begin{abstract}
    We prove the existence of a normal form for a real-analytic Levi-flat hypersurface defined by the vanishing of the real part of a holomorphic function with a Morse-Bott singularity. As a consequence, we recover the Burns-Gong normal form for Levi-flat hypersurfaces with generic Morse singularities and provide a new normal form for a certain class of real analytic quadratic Levi-flat hypersurfaces.
    
    \end{abstract}
\maketitle 

\section{Introduction and Statement of the main result}
\par 
A central result in Morse theory \cite{Morse} is \textit{Morse's Lemma}, which establishes a quadratic normal form for smooth functions in the neighborhood of a non-degenerate critical point. Several generalizations of Morse's Lemma exist; for instance, see Palais \cite{Palais63,Palais69} and Feehan \cite[Theorem 4]{feehan2020}. In the case of holomorphic functions, a version of Morse's Lemma can be found in \cite[p. 102]{Fedo}. 
The condition of \textit{non-degenerate critical points} was relaxed by Bott, who introduced the concept of \textit{non-degenerate manifolds of critical points} in \cite{Bott54}, now known as \textit{Morse-Bott functions}. Bott used this definition in his proof of Bott's Periodicity Theorem \cite{Bott59}. Austin and Braam \cite[Section 3]{Austin} have employed Morse-Bott functions in their approach to developing a Morse-Bott theory for equivariant cohomology. In \cite{Scardua1,Scardua2}, Sc\'ardua and Seade studied codimension-one foliations on closed, oriented manifolds whose singularities are locally defined by Morse-Bott functions. 
A foundational result in this context is the \textit{Morse-Bott Lemma} (see Lemma \ref{morsebott}, \cite[Lemma 3.8]{Petro} or \cite[Corollary 2.15]{feehan2020}). 
\par In this paper, we study \textit{singular real analytic Levi-flat hypersurfaces} and aim to establish a normal form analogous to Morse-Bott Lemma for holomorphic functions. For Levi-flat hypersurfaces, Burns-Gong \cite[Theorem 1.1]{burns} proved the following:  
let $M$  be a germ of real analytic Levi-flat hypersurface at $0\in\mathbb{C}^{n}$, $n\geq 2$, defined by
\[\mathcal{R}e(z_{1}^{2}+\ldots+z_{n}^{2})+H(z,\bar{z})=0\]
with $H(z,\bar{z})=O(|z|^{3})$, $H(z,\bar{z})=\overline{H}(\bar{z},z)$. Then there exists a holomorphic coordinate system such that $$M=\{\mathcal{R}e(x_{1}^{2}+\ldots+x_{n}^{2})=0\}.$$
This result serves as a Morse's Lemma for Levi-flat hypersurfaces, and it is a normal form in the case of a generic (Morse) singularity. Generalizations can be found in \cite{normal}, \cite{singular}, and \cite{marra}.  Recently, new normal forms of Levi-flat hypersurfaces with singularities on a boundary manifold have been obtained in \cite{boundary}. 
\par By combining techniques from holomorphic foliations developed by Cerveau and Lins Neto \cite{alcides} and the Morse-Bott Lemma, we will prove the following theorem.  

\begin{maintheorem}\label{main_theorem}
Let $M=\{ F=0\}$ be a germ of a real analytic Levi-flat hypersurface at $(\mathbb{C}^n,0)$, $n\geq 2$ such that:
\begin{enumerate}
\item $F(z_1,\ldots,z_n)=\mathcal{R}e(z_1^2+\ldots+z^2_{n-c}) + H(z_1,\ldots,z_n,\bar{z}_1,\ldots,\bar{z}_n)$ with $n-c\geq 2$;
\item $H(z,\overline{z})=\overline{H(\overline{z},z)}$, $\frac{\partial{H}}{\partial{z}_{j}}(z,\bar{z})=\frac{\partial{H}}{\partial{\bar{z}}_{j}}(z,\bar{z})=0$ for all $n-c+1\leq j\leq n$, and 
$H(z,\overline{z})=O(|z|^{3})$.
\end{enumerate}
Then there exists a germ of biholomorphism $\Phi\in \operatorname{Diff}(\mathbb{C}^n,0)$ such that
\[D\Phi(0)=\begin{pmatrix}
\operatorname{id}_{n-c} & \star \\
0 & \operatorname{id}_{c} 
\end{pmatrix}\]
where $\operatorname{id}_{n-c}\in\operatorname{GL}(n-c,\C)$, $\operatorname{id}_{c}\in\operatorname{GL}(c,\C)$, and 
\[\Phi^{-1}(M)=\{(x_1,\ldots,x_n)\in(\C^n,0):\mathcal{R}e(x_1^2+\ldots+x_{n-c}^2)=0\}.\]
\end{maintheorem}
\par When $c=0$, Theorem \ref{main_theorem} recovers the result of Burns-Gong \cite[Theorem 1.1]{burns}. 
When $n=3$ and $c=1$, it corresponds to \cite[Theorem 1.3]{line}. Moreover, this theorem yields a new normal form for the real analytic Levi-flat quadratic of type $\mathcal{Q}_{0,2(n-c)}$ given by Burns-Gong \cite[Table 2.1]{burns}.  
\par This paper is organized as follows. In Section \ref{definitions}, we recall some definitions and known results about Levi-flat hypersurfaces and holomorphic foliations. In Section \ref{Morse}, we establish the Morse-Bott Lemma. Finally, Section \ref{theorem1} is devoted to the proof of Theorem \ref{main_theorem}.   
\par The following notations will be used in this paper:
\begin{enumerate}
\item $\mathcal{O}_n$: the ring of germs of holomorphic functions at $0\in\C^n$.
\item $\mathcal{M}_n=\{ f\in\mathcal{O}_n : f(0)=0\}$, the maximal ideal of $\mathcal{O}_n$.
\item $\mathcal{A}_n$: the ring of germs at $0\in\C^n$ of \textit{complex} valued real-analytic functions.
\item $\mathcal{A}_{n\mathbb{R}}$: the ring o germs of \textit{real} valued analytic functions. Note that $f\in\mathcal{A}_n\cap \mathcal{A}_{n\mathbb{R}}\iff f=\overline{f}$.
\item $\dif(\C^n,0)$: the group of germs of biholomorphisms $f:(\C^n,0)\rightarrow (\C^n,0)$ at $0\in\C^n$ with the operation of composition.
\end{enumerate}
\section{Singular Levi-flat hypersurfaces}\label{definitions}
\par Let $M\subset\C^n$, $n\geq 2$, be a germ at the origin of a real-analytic irreducible hypersurface of real codimension one. We may assume $M=\{F(z)=0\}$, with $F\in \mathcal{A}_{n\mathbb{R}}$. The \textit{singular set} of $M$ is given by \[\sing(M)=\{ F(z)=0\}\cap\{dF(z)=0\},\] where $d$ is the usual real differential operator. The \textit{regular part} is defined by $M^*:=M\setminus \sing(M)$. On $M^{*}$, we consider the distribution of complex hyperplanes $L$ given by
\[
L_p:=ker(\partial F(p))\subset T_pM^*=ker(dF(p)),\quad\text{for}\,\,\, p\in M^*.\]
This distribution is called the \textit{Levi distribution} on $M^*$. If $L$ is integrable in the sense of Frobenius, then $M$ is called \textit{Levi-flat}. In this case, $M^*$ is foliated by a real-analytic codimension one foliation $\mathcal{L}$, called the \textit{Levi foliation}. Each leaf of $\mathcal{L}$ is a codimension-one holomorphic submanifold immersed in $M^*$. The Levi distribution can be defined by the real-analytic 1-form $\eta = i(\partial F-\bar{\partial} F)$, the \textit{Levi form} of $F$. The integrability condition is equivalent to the condition
$$
(\partial F-\bpartial F)\wedge \partial\bpartial F\vert_{M^*}=0
$$
or, using the fact that $\partial F +\bpartial F = dF$, is equivalent to
\[
\partial F (p)\wedge \bpartial F(p) \wedge \partial \bpartial F (p) = 0,\,\,\,\,\,\forall\quad p\in M.\]
\par If $\sing(M)=\emptyset$, then $M$ is called \textit{smooth}. In this case, according to E. Cartan \cite[Th\'eor\`eme IV]{Cartan}, around each point $q\in M$ one may find suitable holomorphic coordinates $(z_1,...,z_n)$ of $\C^n$ such that $M$ is locally given by 
\[\{\mathcal{R}e(z_n)=0\}.\] This is the \textit{local normal form} for a smooth real-analytic Levi-flat hypersurface $M$. 
In order to build singular Levi-flat hypersurfaces that are irreducible, we recall the following result from \cite[Lemma 2.2]{alcides} that guaranteed 
the irreducibility of the real-analytic functions. 
\begin{lemma}\label{lema1} Let $f\in\mathcal{M}_n$, $f\neq 0$, and suppose $f$ is not a power in $\mathcal{O}_n$. Then $\mathcal{I}m(f)$ and $\mathcal{R}e(f)$ are irreducible in $\mathcal{A}_{n\mathbb{R}}$. 
\end{lemma} 
\subsection{Complexification of a Levi-flat hypersurface}
Let $F\in\mathcal{A}_{n}$. Its Taylor series at $0\in\mathbb{C}^{n}$ can be written as
\begin{equation}\label{equa-uno}
F(z)=\sum_{\mu,\nu}F_{\mu\nu}z^{\mu}\bar{z}^{\nu},
\end{equation}
 where
$F_{\mu \nu}\in\mathbb{C}$, $\mu=(\mu_{1},\ldots,\mu_{n})$,
$\nu=(\nu_{1},\ldots,\nu_{n})$, $z^{\mu}=z_{1}^{\mu_{1}}\ldots z_{n}^{\mu_{n}}$,
$\bar{z}^{\nu}=\bar{z}_{1}^{\nu_{1}}\ldots \bar{z}_{n}^{\nu_{n}}$. When $F\in\mathcal{A}_{n\mathbb{R}}$, the coefficients $F_{\mu \nu}$ satisfy
$$\bar{F}_{\mu \nu}=F_{\nu \mu}.$$
 The complexification $F_{\mathbb{C}}\in\mathcal{O}_{2n}$ of $F$ is defined by the series
\begin{equation}\label{equa-dos}
F_{\mathbb{C}}(z,w)=\sum_{\mu,\nu}F_{\mu\nu}z^{\mu}w^{\nu}.
\end{equation}
If the series in (\ref{equa-uno}) converges in the polydisc $D^{n}_{r}=\{z\in\mathbb{C}^{n}:|z_{j}|<r\}$ then the
series in (\ref{equa-dos}) converges in the polydisc $D^{2n}_{r}$. Moreover, $F(z)=F_{\mathbb{C}}(z,\bar{z})$ for all
$z\in D^{n}_{r}$.
\par Let $M=\{F=0\}$ be a Levi-flat hypersurface, where $F\in\mathcal{A}_{n\mathbb{R}}$. The complexification $\eta_{\mathbb{C}}$ of its Levi 1-form $\eta=i(\partial{F}-\bar{\partial}F)$ can be written as
$$\eta_{\mathbb{C}}=i(\partial_{z}F_{\mathbb{C}}-\partial_{w}F_{\mathbb{C}})=i
\sum_{\mu,\nu}(F_{\mu\nu}w^{\nu}d(z^{\mu})-F_{\mu\nu}z^{\mu}d(w^{\nu})).$$
\par The complexification $M_{\mathbb{C}}$ of $M$ is
defined as $M_{\mathbb{C}}=\{F_{\mathbb{C}}=0\}$ and its smooth part is
$M^{*}_{\mathbb{C}}=M_{\mathbb{C}}\backslash\{dF_{\mathbb{C}}=0\}$. Clearly $M_{\mathbb{C}}$ defines a complex subvariety of dimension $2n-1$.
The integrability condition of $\eta=i(\partial{F}-\bar{\partial}F)|_{M^{*}}$ implies that $\eta_{\mathbb{C}}|_{M^{*}_{\mathbb{C}}}$ is integrable. Therefore,
$\eta_{\mathbb{C}}|_{M^{*}_{\mathbb{C}}}=0$ defines a holomorphic foliation $\mathcal{L}_{\mathbb{C}}$ on $M^{*}_{\mathbb{C}}$ that will be called the complexification of $\mathcal{L}$.
\begin{remark}\label{alpha-j2}
 Let $\eta=i(\partial{F}-\bar{\partial}F)$ and $\eta_{\mathbb{C}}$ be as above. Then $\eta|_{M^{*}}$ and $\eta_{\mathbb{C}}|_{M^{*}_{\mathbb{C}}}$ define $\mathcal{L}$ and $\mathcal{L}_{\mathbb{C}}$, respectively.
 If we define
$\alpha=\sum_{j=1}^{n}\frac{\partial{F_{\mathbb{C}}}}{\partial{z}_{j}}dz_{j}$ and $\beta=\sum_{j=1}^{n}\frac{\partial{F_{\mathbb{C}}}}{\partial{w}_{j}}dw_{j}$, then $dF_{\mathbb{C}}=\alpha+\beta$
and $\eta_{\mathbb{C}}=i(\alpha-\beta)$, so that
\begin{eqnarray*}
\eta_{\mathbb{C}}|_{M^{*}_{\mathbb{C}}}=
2i\alpha|_{M^{*}_{\mathbb{C}}}=-2i\beta|_{M^{*}_{\mathbb{C}}}.
\end{eqnarray*}
In particular, both $\alpha|_{M^{*}_{\mathbb{C}}}$ and $\beta|_{M^{*}_{\mathbb{C}}}$
define the foliation $\mathcal{L}_{\mathbb{C}}$.
\end{remark}
\subsection{Holomorphic foliations and Levi-flat hypersurfaces}
This section recalls key results regarding Levi-flat hypersurfaces invariant by holomorphic foliations.
\begin{definition}\label{Levi-tangente}
 Let $\mathcal{F}$ and $M=\{F=0\}$ be germs at $(\mathbb{C}^{n},0)$, $n\geq 2$, of a codimension-one
singular holomorphic foliation and of a real Levi-flat hypersurface, respectively. We
say that $\mathcal{F}$ and $M$ are \textit{tangent} if the leaves of the Levi foliation $\mathcal{L}$ on $M$ are also
leaves of $\mathcal{F}$.
\end{definition}
\par We recall that a germ of holomorphic function $h$ is called a holomorphic first integral
for a germ of codimension-one holomorphic foliation $\F$ if its zeros set is contained
in $\sing(\F)$ and its level hypersurfaces contain the leaves of $\F$.
\par The algebraic dimension of $\sing(M)$ is the complex dimension of the singular set of $M_{\mathbb{C}}$. We will use the following result of \cite[Theorem 2]{alcides}, which essentially assures that if the singularities of $M$ are sufficiently small (in the algebraic sense) then $M$ is given by the zeros of the real part of a holomorphic function.
\begin{theorem}[Cerveau-Lins Neto \cite{alcides}]\label{alcides-theorem}
 Let $M=\{F=0\}$ be a germ of an irreducible real analytic Levi-flat hypersurface at $0\in\mathbb{C}^{n}$, $n\geq{2}$, with Levi 1-form $\eta=i(\partial{F}-\bar{\partial} F)$. Assume that the algebraic dimension of $\sing(M)$ is less than or equal to $2n-4$. Then there exists an unique germ at $0\in\mathbb{C}^{n}$ of holomorphic codimension-one foliation $\mathcal{F}_{M}$ tangent to $M$, if one of the following conditions is fulfilled:
\begin{enumerate}
\item[(a)] $n\geq 3$ and $\cod_{M_{\mathbb{C}}^{*}}(\sing(\eta_{\mathbb{C}}|_{M_{\mathbb{C}}^{*}}))\geq 3$.
\item[(b)]  $n\geq 2$, $\cod_{M_{\mathbb{C}}^{*}}(\sing(\eta_{\mathbb{C}}|_{M_{\mathbb{C}}^{*}}))\geq 2$
and $\mathcal{L}_{\mathbb{C}}$ has a non-constant holomorphic first integral.
\end{enumerate}
Moreover, in both cases the foliation $\mathcal{F}_{M}$ has a non-constant holomorphic first integral $f$ such that $M=(\mathcal{R}e(f)=0)$.
\end{theorem}

\subsection{Holonomy and holomorphic first integrals}

Let us now consider a specific situation involving blow-ups and holonomy. Let $\pi:\tilde{\C}^{2n} \rightarrow \C^{2n}$ be the blow-up along a complex submanifold $C\subsetneq\C^{2n}$, and let $E$ denote the exceptional divisor. Denote by $\tilde{M}_\C:=\overline{\pi^{-1}(M_\C \setminus \{C\})} \subset \tilde{\C}^{2n}$ the strict transform of $M_\C$ via $\pi$ and by $\tilde{\mathcal{F}}:=\pi^*(\mathcal{L}_\C)$ the strict transform foliation on $\tilde{M}_\C$ induced by $\mathcal{L}_\C$.
Assume that $\tilde{M}_\C$ is smooth and that the intersection $\tilde{C}=\tilde{M}_\C \cap E$ is invariant by $\tilde{\mathcal{F}}$. Set $S=\tilde{C}\setminus \sing(\tilde{\mathcal{F}})$. Then $S$ is a smooth leaf of $\tilde{\mathcal{F}}$. Fix a point $p_0 \in S$ and a transverse section $\Sigma$ through $p_0$. Let $G \subset \text{Diff} (\Sigma, p_0)$ be the holonomy group of the leaf $S$ of $\tilde{\mathcal{F}}$. Since $\dim_{\C} \Sigma = 1$, we can identify $G \subset \text{Diff} (\C, 0)$. We state the following key lemma.

\begin{lemma}\label{lemaarturo} 
In the above setting, suppose that the following conditions hold:
\begin{enumerate}
\item  For every $p \in S\setminus \sing (\tilde{\mathcal{F}})$, the leaf $L_p$ of $\tilde{\mathcal{F}}$ through $p$ is closed in $S$.
\item $g'(0)$ is a primitive root of unity for all $g\in G \setminus \{\text{id}\}$. 
\end{enumerate}
Then $\mathcal{L}_\C$ admits a non-constant holomorphic first integral.
\end{lemma} 
The proof of this lemma follows the ideas in \cite[Lemma 4.1]{line} and \cite[Theorem 5.1]{singular}.
\section{Morse-Bott Lemma}\label{Morse}
Let $U\subset\C^n$, with $n\geq 2$, be an open subset, and let $f:U\subset\C^n\to\C$ be a holomorphic function. Suppose that the critical set $\operatorname{Crit} f=\{x\in U: f'(x)=0\}$ is a connected holomorphic submanifold. We say that 
$f$ is \textit{Morse-Bott} at a point $x_0\in\operatorname{Crit}f$ if $T_{x_0}\operatorname{Crit}f=\operatorname{Ker} f''(x_0)$, where $f''(x_0)$ is the complex Hessian of $f$ at $x_0$. If $f$ is Morse-Bott at every point $x_0\in\operatorname{Crit}f$, we say that $f$ is \textit{Morse-Bott along} $\operatorname{Crit}f$ or a \textit{Morse-Bott function}.  
\par We now state the holomorphic Morse-Bott Lemma; see for instance \cite[Corollary 2.15]{feehan2020} or \cite[Lemma 3.8]{Petro}.
\begin{lemma}\label{morsebott}  
Let $n\geq2$ be an integer, $U\subset\C^n$ be an open neighborhood of the origin, and $f:U\to\C$ be a holomorphic function such that $f(0)=0$ and $f'(0)=0$. Assume that $\text{Crit}(f)$ is a complex manifold of $U$ with complex tangent space $T_0 \operatorname{Crit} f=\operatorname{Ker} f''(0)$ of dimension $c\geq 0$ at the origin. Then, after possibly shrinking $U$, there are an open neighborhood $V\subset\C^n$ of the origin and a biholomorphism,
\[
V\ni(x_1,\ldots,x_n)\mapsto(z_1,\ldots,z_n)=\Phi(x_1,\ldots,x_n)\in\C^n
\]
onto an open neighborhood of the origin in $\C^n$ such that
\[\Phi^{-1}(U\cap\operatorname{Crit} f)=V\cap(\C^{c}\times\{0\})\subset\C^{c}\times\C^{n-c}\]
with $\Phi(0)=0$ and 
\[D\Phi(0)=\begin{pmatrix}
\operatorname{id}_{n-c} & \star \\
0 & \operatorname{id}_{c} 
\end{pmatrix}\]
where $\operatorname{id}_{n-c}\in\operatorname{GL}(n-c,\C)$, $\operatorname{id}_{c}\in\operatorname{GL}(c,\C)$ and 
\[f(\Phi(x_1,\ldots,x_n))=x_1^2+\ldots+x_{n-c}^2,\,\,\,\,\,\text{for all}\,\,\,x=(x_1,\ldots,x_n)\in U.\]
\end{lemma}
This lemma plays a key role in the proof of Theorem \ref{main_theorem}, providing the desired normal form in a neighborhood of a Morse-Bott singularity.

\section{Proof of Theorem \ref{main_theorem}}\label{theorem1}
Let $M=\{ F=0\}\subset\C^n$, with $n\geq 2$, be a germ of a real analytic Levi-flat hypersurface at the origin, where
\[F(z_1,\ldots,z_n)=\mathcal{R}e(z_1^2+\ldots+z^2_{n-c}) + H(z_1,\ldots,z_n,\bar{z}_1,\ldots,\bar{z}_n)\,\,\,\,\text{with}\,\,\,\,\,n-c\geq 2,\]
and assume that $H(z,\overline{z})$ satisfies the hypothesis of Theorem \ref{main_theorem}.
\par The proof is based on applying Theorem \ref{alcides-theorem} to obtain a germ $f \in \mathcal{O}_n$ such that the foliation $\F_M$ defined by $\omega=d f$ is tangent to $M$ and 
\begin{equation}\label{eq_1}
M= \{(z_1,\ldots,z_n)\in(\C^n,0): \mathcal{R}e(f(z_1,\ldots,z_n))=0\}. 
\end{equation}
Let us assume the existence of such a function germ and proceed to complete the proof. 
The foliation $\mathcal{F}_M$ can be viewed as an extension to a neighborhood of $0 \in \C^n$ of the Levi foliation $\mathcal{L}$ on $M^*$.
Without loss of generality, assume that $f$ is not a power in $\mathcal{O}_n$. Then $\mathcal{R}e(f)$ is irreducible by \cite[Lemma 2.2]{alcides}.
Since $M=\{F=0\}$, the equation (\ref{eq_1}) implies that $\mathcal{R}e(f)=U\cdot F$, where $U\in\mathcal{A}_{n\mathbb{R}}$ and $U(0)\neq 0$. In particular, we get  
\begin{equation*}
f(z_1,\ldots,z_n)=U(0)(z_1^2+\cdots+z^2_{n-c})+\mathcal{O}(|z|^{3}).
\end{equation*}
Now, we can apply Lemma \ref{morsebott} to $f$, which implies the existence of a biholomorphism germ $\Phi\in \operatorname{Diff}(\mathbb{C}^n,0)$ such that
\[D\Phi(0)=\begin{pmatrix}
\operatorname{id}_{n-c} & \star \\
0 & \operatorname{id}_{c} 
\end{pmatrix}\]
where $\operatorname{id}_{n-c}\in\operatorname{GL}(n-c,\C)$, $\operatorname{id}_{c}\in\operatorname{GL}(c,\C)$ and
\[f(\Phi(x_1,\ldots,x_n))=x_1^2+\ldots+x_{n-c}^2,\,\,\,\,\,\text{for all}\,\,\,w=(x_1,\ldots,x_n)\in (\C^n,0).\]
Thus
\begin{eqnarray*}
\Phi^{-1}(M)&=&\{(x_1,\ldots,x_n)\in (\C^n,0):\Phi(x_1,\ldots,x_n)\in M\}\\
&=& \{(x_1,\ldots,x_n)\in (\C^n,0): \mathcal{R}e(f(\Phi(x_1,\ldots,x_n)))=0\}\\
&=& \{(x_1,\ldots,x_n)\in (\C^n,0): \mathcal{R}e(x_1^2+\ldots+x_{n-c}^2)=0\},
 \end{eqnarray*}
 and Theorem \ref{main_theorem} is proved.
\par Next, we will show that $M=\{F(z)=0\}$ satisfies the hypotheses of Theorem \ref{alcides-theorem}. 
We have
\[F(z)=\mathcal{R}e(z_1^2+\ldots+z^2_{n-c}) + H(z,\bar{z}),\]
where $H(z,\overline{z})=\overline{H(\overline{z},z)}$ and 
$H(z,\overline{z})=O(|z|^{3})$. Then, the complexification $F_{\C}$ of $F$ is given by 
\[
F_\C(z,w)=\frac{1}{2}(z_1^2+\ldots+z^2_{n-c})+\frac{1}{2}(w_1^2+\ldots+w^2_{n-c})+H_\C(z,w),\]
where $H_\C(z,\bar{z})=H(z,\bar{z})$. Therefore, 
$M_\C=\{(z,w)\in(\C^{2n},0): F_\C(z,w)=0\}$ is a germ of complex analytic subvariety of dimension $2n-1$ at $0\in\C^{2n}$ whose singular set is \[\sing(M_\C)=\{z_1=\ldots=z_{n-c}=w_1=\ldots=w_{n-c}=0\}.\]
Since $\sing(M_\C)$ has dimension $2c$, and $n-c\geq 2$, the algebraic dimension of $\sing(M)$ is $\leq 2n-4$, satisfying the first hypothesis of Theorem \ref{alcides-theorem}. On the other hand, the complexification of $\eta = i(\partial F - \overline{\partial} F)$ is
$\eta_\C = i(\partial_z F_{\C} - {\partial}_w F_{\C})$.
Recall that $\eta\vert_{M^*}$ and $\eta_{\C}\vert_{M^*_\C}$ define $\mathcal{L}$ and $\mathcal{L}_\C$ respectively. Now we compute $\sing(\eta_\C\vert_{M^*_\C})$. Since $\frac{\partial{H}}{\partial{z}_{j}}(z,\bar{z})=\frac{\partial{H}}{\partial{\bar{z}}_{j}}(z,\bar{z})=0$ for all $n-c+1\leq j\leq n$, we can write $dF_\C = \alpha + \beta$ where
$$
\alpha=\sum^n_{j=1} \dfrac{\partial F_\C}{\partial z_j} = \sum_{j=1}^{n-c} z_i \, dz_i + \sum_{j=1}^{n-c} \dfrac{\partial H_\C}{\partial z_j} \, dz_j 
$$
and
$$
\beta=\sum^n_{j=1} \dfrac{\partial F_\C}{\partial w_j} = \sum_{j=1}^{n-c} w_i \, dw_i + \sum_{j=1}^{n-c} \dfrac{\partial H_\C}{\partial w_j} \, dw_j.
$$
Then $\eta_\C = i(\alpha - \beta)$ and so
$$
\eta_\C \vert_{M^*_\C} = (\eta_\C + idF_\C)\vert_{M^*_\C} = 2i\alpha\vert_{M^*_\C} = -2i\beta\vert_{M^*_\C}.
$$
In particular, $\alpha\vert_{M^*_\C}$ and $\beta\vert_{M^*_\C}$ define $\mathcal{L}_\C$. Therefore, $\sing(\eta_\C \vert_{M^*_\C})$ can be split in two parts. 
Let $A_j=\dfrac{\partial H_\C}{\partial z_j}$ and $B_j = \dfrac{\partial H_\C}{\partial w_j}$. Let \[M_1 = \{(z, w) \in M_\C:\quad\frac{\partial F_\C}{\partial w_j}\neq 0 \text{ for some } j = 1, . . . , n\}\] and \[M_2 = \{(z, w) \in M_\C:\quad \frac{\partial F_\C}{\partial z_j}\neq 0 \text{ for some } j = 1, . . . , n\}.\] Note that $M_\C = M_1 \cup M_2$; if we denote by
$$
X_1:= M_1 \cap \left\{ z_1 +A_1 = ... = z_{n-c} +A_{n-c}=0 \right\}
$$
and
$$
X_2:= M_2 \cap \left\{ w_1 +B_1 = ... = w_{n-c} +B_{n-c}=0 \right\}
$$
then $\sing(\eta_\C\vert_{M^*_\C}) = X_1 \cup X_2$ and clearly
$cod_{M^*_\C}\sing(\eta_\C\vert_{M^*_\C})=n-c$. 
If $n-c\geq 3$, we can directly apply Theorem \ref{alcides-theorem} $(a)$ and the proof ends. 
\par In the case $n-c = 2$, we are going to prove that $\mathcal{L}_{\C}$ has a non-constant holomorphic first integral.  
Let $F(z)=\mathcal{R}e(z_1^2+z^2_{2}) + H(z,\bar{z})$,
where $H(z,\overline{z})=\overline{H(\overline{z},z)}$ and 
$H(z,\overline{z})=O(|z|^{3})$. Then, the complexification $F_{\C}$ of $F$ is given by 
\[
F_\C(z,w)=\frac{1}{2}(z_1^2+z^2_{2})+\frac{1}{2}(w_1^2+w^2_{2})+H_\C(z,w),\]
where $H_\C(z,\bar{z})=H(z,\bar{z})$.
Consider the blow-up $\pi:\tilde{\C}^{2n}\to\C^{2n}$ along
$$
C = \{ z_1 = z_2 = w_1 = w_2 = 0 \} \simeq \C^{2(n-2)}
$$
with exceptional divisor $E$.
Let $\tilde{M}_\C := \overline{\pi^{-1} ( M_\C \setminus \{C\} ) } \subset \tilde{\C}^{2n}$ be the strict transform of $M_\C$ by $\pi$ and $\tilde{\mathcal{F}} := \pi^*(\mathcal{L}_\C)$ be the strict transform foliation on $\tilde{M}_\C$.

Consider, for instance, the chart $(U_1, (t,s)=(t_1, t_2,\ldots,t_n, s_1,s_2,\ldots,s_n))$ of $\tilde{\C}^{2n}$. For convenience, we will denote $t_1 = u$. In this chart,
\[
\pi(u,t_2,\ldots,t_n,s_1,\ldots,s_n) = (u,u\, t_2, t_3,\ldots,t_n, u\, s_1, u \, s_2, s_3,\ldots,s_n).\]
We have
\[
\tilde{M}_\C \cap U_1 = \{ 1 + t_2^2+s_1^2+s_2^2 + u\cdot H_1(t,s) = 0\},\]
where
$$
H_1(t,s) = \frac{H(t,s)}{u}
$$
and from this
$$
E \cap \tilde{M}_\C \cap U_1 = \{ 1 + t_2^2+s_1^2+s_2^2 = u = 0\}
$$
we get that these complex subvarieties are smooth. The foliation $\tilde{\mathcal{F}}$ on $U_1$ is defined by $\tilde{\beta}\vert_{\tilde{M}_\C \cap U_1}=0$, where

\begin{equation}\label{betatil}
\tilde{\beta} = (s_1^2+s_2^2) \, du + us_1\, ds_1 + u s_2 \, ds_2 + u\, \tilde{\theta},
\end{equation}
and
$$
\tilde{\theta} = \dfrac{\pi^* \left( \displaystyle\sum_{j=1}^{2} \frac{1}{2} B_j \, dw_j \right)}{u^2}.
$$
Note that the exceptional divisor $E$ is invariant by $\tilde{\mathcal{F}}$ and the intersection with $\sing(\tilde{\mathcal{F}})$ is
$$
\sing(\tilde{\mathcal{F}}\vert_{U_1}) \cap E = \{ u = s_1^2 + s_2^2 = 1+t_2^2\}.
$$
In particular, $S=(E\cap \tilde{M}_\C)\setminus \sing(\tilde{\mathcal{F}})$ is a leaf of $\tilde{\mathcal{F}}$. We calculate the generators of the holonomy group $G$ of the leaf $S$. Due to the symmetry of the variables in the definition of $\tilde{M}_\C$, working in the chart $U_1$ is sufficient.

Pick $p_0 = (0,\ldots,0,i,0,\ldots,0) \in S \cap U_1$ (that is, $s_1=i$) and the transverse
$ \Sigma = \{ (\lambda,0,\ldots,0,i,0,\ldots,0): \lambda \in \mathbb{C} \}$ parameterized by $\lambda$  at $p_0$. Let $\rho_1, \rho_2$ be $2^{nd}$-primitive roots of $-1$. We have 
$$
\begin{array}{rcl}
\sing(\tilde{\mathcal{F}}\vert_{U_1})\cap E &=& \{ u = s_1+is_2 = t_2-\rho_1=0 \}\cup 
  \{ u = s_1-is_2 = t_2-\rho_1=0 \} \cup\\
& &  \{ u = s_1+is_2 = t_2-\rho_2=0 \} \cup
  \{ u = s_1-is_2 = t_2-\rho_2=0 \}.
\end{array}
$$
The fundamental group $\pi_1(S,p_0)$ can be written in terms of generators as
$$
\pi_1(S,p_0)= \langle \gamma_j, \delta_j \rangle_{j=1,2},
$$
where, for each $j=1,2$, $\gamma_j$ are loops around the connected component of the singular set
$$
\{ u = s_1+is_2 = t_2-\rho_j=0 \}
$$
and $\delta_j$ are loops around the connected component of the singular set
$$
\{ u = s_1-is_2 = t_2-\rho_j=0 \}.
$$
Therefore, $G=\langle f_j, g_j\rangle_{j=1,2}$, where $f_j$ and $g_j$ correspond to $[\gamma_j]$ and $[\delta_j]$, respectively. We get from \eqref{betatil} that $f'_j(0)=\text{e}^{-2\pi i}$ and $g'_j(0)=\text{e}^{-2\pi i}$, for $j=1,2$. 
Finally, Lemma \ref{lemaarturo} implies that $\mathcal{L}_\C$ has a non-constant holomorphic first integral, completing the proof.

\begin{exemplo} Consider $$M=\left\{ (z_1,z_2,z_3,z_4)\in (\mathbb{C}^4,0) : \mathcal{R}e(z_1^2+z_2^2)+z_1\bar{z}_1z_2\bar{z}_2=0 \right\}.$$ 
Setting $H(z,\bar{z})=z_1\bar{z}_1z_2\bar{z}_2 \in \mathcal{O}(|z|^4)$, one easily sees that $\frac{\partial H}{\partial{z}_i}=\frac{\partial H}{\partial \bar{z}_i} = 0$ for $i=3,4$. Theorem 1 asserts that there exists a germ of biholomorphism $\Phi\in \operatorname{Diff}(\mathbb{C}^4,0)$ such that $$\Phi^{-1} (M) = \left\{ (x_1,x_2,x_3,x_4) \in (\mathbb{C}^4,0): \mathcal{R}e(x_1^2+x_2^2)=0\right\}.$$
\end{exemplo}

\vspace{1cm}

\noindent{\bf Data Availability Statement:} 
Data sharing is not applicable to this article as no data sets were generated or analyzed during the current study.
\vspace{1cm}

\noindent{\bf Declarations
Conflict of Interest:} The authors declare that they have no conflict of interest.

\end{document}